\documentclass{article}    
\usepackage{amstext}

\title{Extending the idea of compressed algebra to arbitrary socle-vectors}  
\author{Fabrizio Zanello\\
Department of Mathematics and Statistics\\
Queen's University\\
Kingston, Ontario\\
K7L 3N6 Canada\\
E-mail: fabrizio@mast.queensu.ca\\
\\}         

\def\qed{$\rlap{$\sqcap$}\sqcup$}
\usepackage{latexsym}

\begin{document}           

{\large

\maketitle                 

{\ }\\

\section{Introduction}

In this paper we will study artinian quotients $A=R/I$ of the polynomial ring $R=k[x_1,...,x_r]$, where $k$ is a field of characteristic zero, the $x_i$'s all have degree 1 and $I$ is a homogeneous ideal of $R$. These rings are often referred to as standard graded artinian algebras.\\
Before explaining the main results of this work, we establish some of the notation we will use: the {\it h-vector} of $A$ is $h(A)=h=(h_0,...,h_e)$, where $h_i=\dim_k A_i$ and $e$ is the last index such that $\dim_k A_e>0$. Since we may suppose that $I$ does not contain non-zero forms of degree 1, $r=h_1$ is defined to be the {\it embedding dimension} ({\it emb.dim.}, in brief) of $A$.\\
The {\it socle} of $A$ is the annihilator of the maximal homogeneous ideal $\overline{m}=(\overline{x_1},...,\overline{x_r})\subset A$, namely $soc(A)=\lbrace a\in A {\ } \mid {\ } a\overline{m}=0\rbrace $. Since $soc(A)$ is a homogeneous ideal, we define the {\it socle-vector} of $A$ as $s(A)=s=(s_0,...,s_e)$, where $s_i=\dim_k soc(A)_i$. Note that $s_e=h_e>0$.\\
We will say that an $h$-vector $h$ is {\it admissible for the pair $(r,s)$} if there exists an algebra $A$ with emb.dim.($A$)$=r$, $s(A)=s$ and $h(A)=h$. When the pair $(r,s)$ is clear from the context, we will simply say that $h$ is {\it admissible}.\\
A natural question which arises is the following: what are the admissible $h$-vectors for a given pair $(r,s)$?\\
\\
This problem has been considered in several different guises: e.g. there are several papers which treat the question of determining the possible $h$-vectors of Gorenstein algebras, i.e. finding all the admissible $h$'s which correspond to a fixed $r$ and $s=(0,...,0,1)$. See e.g. Stanley ($[St]$, Thm. 4.2), who gives a complete answer to this question for the case $r\leq 3$ (the case $r=2$ was actually already known to Macaulay, cf. $[Ma]$). See also Migliore and Nagel's paper $[MN]$ for the $h$-vectors of Gorenstein algebras with the Weak Lefschetz Property.\\
More generally, there are many papers which consider the question for the case of level algebras (e.g. see $[BG]$, $[CI]$).\\
The problem of finding {\it all} the admissible $h$-vectors for a given pair $(r,s)$ seems very difficult in general. Iarrobino (cf. $[Ia]$) and Fr\"oberg and Laksov (cf. $[FL]$) considered a more restricted question. More precisely, Iarrobino showed that, putting some natural restrictions on a given pair $(r,s)$, any admissible $h$-vector is bounded from above by a certain maximal $h$, and defined an algebra $A$ with the data $(r,s)$ as {\it compressed} if this maximal $h$ satisfies $h=h(A)$; moreover he proved that, under his hypotheses on $r$ and $s$, there always exists a compressed algebra.\\
This problem of Iarrobino's was taken up again in $[FL]$ by Fr\"oberg and Laksov, who used a different approach.\\
We finally recall the seminal work on compressed algebras, Emsalem and Iarrobino's 1978 article $[EI]$.\\
In this paper we take a more general view. We ask the question: given {\it any} $(r,s)$, is there an $h$ which is maximal among all the admissible $h$-vectors? If such an $h$ exists, we will define as {\it generalized compressed} any algebra with the data $(r,s,h)$ (see Section 2). Naturally, this more general definition coincides with Iarrobino's in the cases satisfying his conditions, and, with our generalized definition, we are enlarging the set of compressed algebras beyond those found in $[Ia]$ and $[FL]$.\\
\\
Let us fix the emb.dim. $r$ and the socle-vector $s=(s_0=0,s_1,...,s_e)$. The two main results of this paper are: Theorem A, an upper-bound $H$ for the $h$-vectors admissible for the pair $(r,s)$, which improves the one given by Fr\"oberg and Laksov in $[FL]$; Theorem B, which asserts that, under certain conditions on $(r,s)$ (less strong than those of $[Ia]$ and $[FL]$), there exist algebras having exactly the upper-bound $H$ we described above.\\
In some cases, however, we will see that the $H$ given by Theorem A is not admissible, and we will supply counter-examples. These counter-examples, moreover, show that the hypotheses of Theorem B, in general, cannot be improved.\\
Here we only mention that, in our forthcoming paper $[Za]$, we will prove that a generalized compressed algebra does {\it not} exist for every pair $(r,s)$, where $r$ is greater than or equal to the {\it minimum embedding dimension} of $s$ (briefly, {\it min.emb.dim.$(s)$}), i.e. the least emb.dim. $r$ such that there exists any algebra $A$ with data $(r,s)$.\\
\\
The results obtained in this paper will be part of the author's Ph.D. dissertation, written at Queen's University (Kingston, Ontario, Canada), under the supervision of Professor A.V. Geramita.
\\
\\
\section{Definitions and preliminary results}

Fix $r$ and $s=(s_0=0,s_1,...,s_e)$; from now on we may suppose, to avoid trivial cases, that $r>1$ and $e>1$.\\
\\
{\bf Definition-Remark 2.1.} Following $[FL]$, define, for $d=0,1,...,e$, the integers $$r_d=N(r,d)-N(r,0)s_d-N(r,1)s_{d+1}-...-N(r,e-d)s_e,$$ where $$N(r,d)=\dim_kR_d={r-1+d \choose d}.$$
It is easy to show (cf. $[FL]$) that $r_0<0$, $r_e\geq 0$ and $r_{d+1}>r_d$ for every $d$.\\
Define $b$, then, as the unique index such that $1\leq b\leq e$, $r_b\geq 0$ and $r_{b-1}<0$.\\
\\
Let $S=k[y_1,...,y_r]$, and consider $S$ as a graded $R$-module where the action of $x_i$ on $S$ is partial differentiation with respect to $y_i$.\\
Recall that, in the theory of Inverse Systems, the $R$-submodule $M$ of $S$ associated to the algebra $R/I$ with data $(r,s)$ is generated by $s_i$ elements of degree $i$, for $i=1,...,e$, and the $h$-vector of $R/I$ is given by the number of linearly independent derivatives in each degree obtained by differentiating the generators of $M$.\\
The number $$N(r,d)-r_d=N(r,0)s_d+N(r,1)s_{d+1}+...+N(r,e-d)s_e$$ is an upper-bound for the number of linearly independent derivatives yielded in degree $d$ by the generators of $M$ and, therefore, is also an upper-bound for the $h$-vector of $R/I$. This is the reason for the introduction of the numbers $r_d$.\\
For a complete introduction to Inverse Systems, we refer the reader to $[Ge]$.\\
\\
{\bf Remark 2.2.} It is easy to see that, for any pair $(r,s)$, we must have $$b\geq e/2.$$ In fact, $$r_b=N(r,b)-N(r,0)s_b-N(r,1)s_{b+1}-...-N(r,e-b)s_e.$$ If $b<e/2$, then $b<e-b$, hence $N(r,e-b)>N(r,b)$. Since $s_e\geq 1$, we get $r_b<0$, a contradiction.\\
\\
{\bf Proposition 2.3} (Fr\"oberg-Laksov). Let $(r,s)$ be as above, $r\geq $ min.emb.dim.$(s)$. Then an upper-bound for the $h$-vectors admissible for the pair $(r,s)$ is given by $$H=(h_0,h_1,...,h_e),$$ where $$h_i=\min \lbrace N(r,i)-r_i,N(r,i)\rbrace$$ for $i=0,1,...,e$.\\
\\
{\bf Proof.} See $[FL]$, Prop. 4, $i$).{\ }{\ }\qed \\
\\
{\bf Remark 2.4.} Fr\"oberg and Laksov have given a direct proof of the proposition; a second proof follows immediately from our comment about Inverse Systems and the numbers $r_d$. The same upper-bound was already supplied by Iarrobino (cf. $[Ia]$) under the natural restriction $s_1=...=s_{b-1}=0$.\\
\\
{\bf Lemma 2.5} (Iarrobino, Fr\"oberg-Laksov). Let $S=k[y_1,...,y_r]$ be the $R$-module defined above, and consider $n$ generic forms $F_1,...,F_n\in S$, respectively of degrees $d_1,...,d_n$. Then, for every integer $c\geq 0$, the subspace of $S_c$ spanned by $R_{d_1-c}F_1,...,R_{d_n-c}F_n$ has dimension (as a $k$-vector space) equal to $$\min \lbrace \dim_k S_c, \dim_k S_{d_1-c}+...+\dim_k S_{d_n-c} \rbrace ,$$ i.e. generic forms have derivatives as independent as they can be.\\
\\
{\bf Proof.} See $[Ia]$, Prop. 3.4 and $[FL]$, Prop. 20. The case $n=1$ was already known to Emsalem and Iarrobino (see $[EI]$).{\ }{\ }\qed \\
\\
From the previous lemma we immediately obtain:\\
\\
{\bf Proposition 2.6} (Iarrobino, Fr\"oberg-Laksov). Let $(r,s)$ be as above, $r\geq $ min.emb.dim.$(s)$. If, moreover, $s_1=...=s_{b-1}=0$, then the upper-bound $H$ yielded by Proposition 2.3 is admissible for the pair $(r,s)$.\\
\\
{\bf Proof.} See $[Ia]$, Thm. II A; $[FL]$, Prop. 4, $iv$) and Thm. 14.{\ }{\ }\qed \\
\\
{\bf Definition 2.7.} 
\begin{itemize}
\item[  $i$)] Fix a pair $(r,s)$ such that $s_1=...=s_{b-1}=0$. Iarrobino (cf. $[Ia]$) defined an algebra as {\it compressed} with respect to this pair $(r,s)$ if its $h$-vector is the upper-bound $H$ of Proposition 2.3.
\item[  $ii$)] Now fix {\it any} pair $(r,s)$. We define an algebra as {\it generalized compressed} with respect to the pair $(r,s)$ if its $h$-vector is the maximal among all the admissible $h$-vectors.
\end{itemize}
{\ }\\
{\bf Remark 2.8.} Proposition 2.6 shows that, for the pairs $(r,s)$ as in Iarrobino's definition (see Definition 2.7, $i$)), there always exists a compressed algebra. However, from Remark 2.2 we see that the restrictions required on the socle-vector $s$ in order to satisfy Iarrobino's conditions are very strong (at least the first half of $s$ must be zero). This is one of the main reasons that lead us to extend the concept and look for generalized compressed algebras.\\
\\
{\bf Example 2.9.} Let $r=3$, $s=(0,0,0,0,0,2,3,1)$. It is easy to see that $r_4=-19<0$ and $r_5=4\geq 0$, whence $b=5$. Since $s_1=...=s_4=0$, Proposition 2.6 says that there is a compressed algebra for this pair $(r,s)$, having $h$-vector $$H=(1,3,6,10,15,17,6,1).$$
The method used by $[Ia]$ and $[FL]$ (suggested by Lemma 2.5) to construct such a compressed algebra is the following: choose one generic form of degree 7 (yielding 3 linearly independent first derivatives and 6 linearly independent second derivatives), three generic forms of degree 6 (yielding 9 first derivatives), and two generic forms of degree 5. Then, by Lemma 2.5, the total number of linearly independent derivatives supplied in degree 4 is $=\min \lbrace 10+18+6, N(3,4)=15\rbrace =15$, whence we obtain that these derivatives span exactly $S_4$; thus we have constructed our $H$.\\
\\
The next two examples illustrate some of the limitations inherent in Propositions 2.3 and 2.6.\\
\\
{\bf Example 2.10.} Let $r=3$, $s=(0,0,0,0,0,3,3,1,2)$. It is easy to check that $b=6$. By Proposition 2.3, the upper-bound for the admissible $h$-vectors for this pair $(r,s)$ is $$H=(1,3,6,10,15,21,18,7,2).$$ Is $H$ admissible? Proposition 2.6 gives no information, since $s_{b-1}=s_5=3>0$. We will see later (as a consequence of our Theorem B) that the answer is positive. We will also show how to construct a (generalized compressed) algebra with $h$-vector $H$.\\
\\
{\bf Example 2.11.} Let $r=3$, $s=(0,0,1,2,1,1,1,1)$. Here $b=5$ and, by Proposition 2.3, the upper-bound for the admissible $h$-vectors for this pair $(r,s)$ is $$H=(1,3,6,10,15,10,4,1).$$ Is $H$ admissible? Proposition 2.6 gives no information, since $s_{b-1}=s_4=1>0$. We will see, by our Theorem A, that in this case the answer is negative, i.e. $H$ is not admissible, since there is an upper-bound sharper than $H$ for this pair $(r,s)$.\\
\\
{\bf Definition-Remark 2.12.} Let $n$ and $i$ be positive integers. The {\it i-binomial expansion of n} is $$n_{(i)}={n_i\choose i}+{n_{i-1}\choose i-1}+...+{n_j\choose j},$$ where $n_i>n_{i-1}>...>n_j\geq j\geq 1$.\\
Under these hypotheses, the $i$-binomial expansion of $n$ is unique.\\
Following $[BG]$, define, for any integer $a$, $$(n_{(i)})_{a}^{a}={n_i+a\choose i+a}+{n_{i-1}+a\choose i-1+a}+...+{n_j+a\choose j+a}.$$
\\
A well-known result of Macaulay is:\\
\\
{\bf Theorem 2.13} (Macaulay). Let $h=(h_i)_{i\geq 0}$ be a sequence of non-negative integers, such that $h_0=1$, $h_1=r$ and $h_i=0$ for $i>e$. Then $h$ is the $h$-vector of some standard graded artinian algebra if and only if, for every $1\leq d\leq e-1$, $$h_{d+1}\leq ((h_d)_{(d)})_{+1}^{+1}.$$
\\
{\bf Proof.} See $[St]$.{\ }{\ }\qed \\
\\
{\bf Remark 2.14.} This result actually holds, with an analogous statement, for any standard graded algebra, not necessarily artinian.\\
\\
{\bf Lemma 2.15} (Bigatti-Geramita). Let $a,b$ be positive integers, $b>1$. Then the smallest integer $s$ such that $a\leq (s_{(b-1)})_{+1}^{+1}$ is $$s=(a_{(b)})_{-1}^{-1}.$$
\\
{\bf Proof.} See $[BG]$, Lemma 3.3.{\ }{\ }\qed \\
\\
{\bf Remark 2.16.} This result yields a lower-bound for the $i$-th entry of an $h$-vector, once the $(i+1)$-st entry is known. In terms of Inverse Systems, it supplies a lower-bound for the number of linearly independent first derivatives of any given set of linearly independent forms of degree $i+1$.\\
\\
Now we state two fundamental results of Iarrobino about sums of powers of linear forms. They will play a key role in the proof of our Theorem B.\\
\\
{\bf Theorem 2.17} (Iarrobino). Let $F=\sum_{t=1}^m L_t^d$ be a form of degree $d$ in $S=k[y_1,...,y_r]$, where the $L_t=\sum_{k=1}^rb_{tk}y_k$ are linear forms, and let $I\subset R$, $I=Ann(F)$. Then there exists a non-empty open subset $U$ of $k^{rm}$ such that, for any choice of the coefficients $b_{tk}$ from $U$, the Gorenstein artinian algebras $R/I$ all have the same $h$-vector, denoted by: $$h(m,d)=(1,h_1(m,d),...,h_d(m,d)=1)=(1,h_1,...,h_d=1),$$ where $$h_s=\min \lbrace m,\dim_kR_s,\dim_kR_{d-s}\rbrace .$$
\\
{\bf Proof.} See $[Ia]$, Prop. 4.7.{\ }{\ }\qed \\
\\
{\bf Theorem 2.18} (Iarrobino). For $i=1,...,n$, let $F_i=\sum_{t=1}^{m_i}L_{it}^{d_i}$ be sums of powers of generic linear forms as above, with $d_1\leq ...\leq d_n$. Then the algebra $A=R/I$, where $I=Ann(\sum_{i=1}^n(F_i))$, has $h$-vector $$h=(1,h_1,...,h_{d_n}),$$ where $$h_s=\min \lbrace \dim_kR_s, \sum_{i=1}^nh_s(m_i,d_i) \rbrace ,$$ provided that $\sum_{i=1}^nh_s(m_i,d_i)\leq \dim_kR_s$ for each $s=d_1,d_2,...,d_n$.\\
\\
{\bf Proof.} See $[Ia]$, Thm. 4.8 B.{\ }{\ }\qed \\
\\
\\
\section{The main results}

We are now ready for the two main results of this paper. The first, as we already mentioned, is an upper-bound $H$ for all the $h$-vectors admissible for given emb.dim. $r$ and socle-vector $s=(s_0=0,s_1,...,s_e)$:\\
\\
{\bf Theorem A.} Let $(r,s)$ be as above, $r\geq $ min.emb.dim.$(s)$. Then an upper-bound $H$ for the $h$-vectors admissible for the pair $(r,s)$ is given by $$H=(h_0,h_1,...,h_e),$$ where $h_0=1$, $h_1=r$ and, inductively, for $2\leq i\leq e$, $$h_i=\min \lbrace ((h_{i-1}-s_{i-1})_{(i-1)})_{+1}^{+1}, N(r,i)-r_i\rbrace .$$
\\
{\bf Proof.} $h_0=1$ and $h_1=r$ is obvious. By induction, let the theorem hold up to some $j$, $1\leq j\leq e-1$. Then, by Proposition 2.3, $h_{j+1}\leq N(r,j+1)-r_{j+1}$, and, using Inverse Systems, by Theorem 2.13, the largest number of linearly independent forms of degree $j+1$ having $h_j-s_j$ first derivatives is $((h_j-s_j)_{(j)})_{+1}^{+1}$. This concludes the induction and the proof of the theorem.{\ }{\ }\qed \\
\\
{\bf Remark 3.1.} Note that, once $h_i=N(r,i)-r_i$ (i.e. $N(r,i)-r_i\leq ((h_{i-1}-s_{i-1})_{(i-1)})_{+1}^{+1}$), then also $h_{i+1}$ has the same form, i.e. $h_{i+1}=N(r,i+1)-r_{i+1}$: in fact $$((N(r,i)-r_i-s_i)_{(i)})_{+1}^{+1}\geq N(r,i)-r_i-s_i=$$ $$N(r,1)s_{i+1}+N(r,2)s_{i+2}+...+N(r,e-i)s_e\geq $$ $$N(r,0)s_{i+1}+N(r,1)s_{i+2}+...+N(r,e-i-1)s_e=N(r,i+1)-r_{i+1}.$$\\
\\
{\bf Remark 3.2.} In general, our upper-bound is sharper than that supplied by Fr\"oberg-Laksov. In fact, it is easy to see, by induction, that $$((h_{i-1}-s_{i-1})_{(i-1)})_{+1}^{+1}\leq N(r,i),$$ since, at each step, $h_{i-1}\leq N(r,i-1)$, $s_{i-1}\geq 0$ and the function $(n_{(t)})_{+1}^{+1}$ increases with $n$.\\
For instance, let $r=3$ and $s=(0,0,1,2,1,1,1,1)$, as in Example 2.11. Then, as we already saw, the upper-bound given by Proposition 2.3 is $$H=(1,3,6,10,15,10,4,1);$$ instead Theorem A yields the sharper $$H=(1,3,6,7,6,6,4,1).$$
\\
More precisely, we have the following:\\
\\
{\bf Proposition 3.3.} The upper-bounds $H$ yielded by Proposition 2.3 and Theorem A are the same if and only if $s_0=s_1=...=s_{b-2}=0$ and $$s_{b-1}\leq N(r,b-1)-((N(r,b)-r_b)_{(b)})_{-1}^{-1}.$$
Otherwise, Theorem A yields a sharper $H$.\\
\\
{\bf Proof.} By Remark 3.2, it remains only to show the first assertion. If the two vectors $H$ are the same, then, for $i=2,...,b-1$, $$((h_{i-1}-s_{i-1})_{(i-1)})_{+1}^{+1}=N(r,i).$$ Therefore, by induction and the properties of the binomial expansion, we have at once $s_1=...=s_{b-2}=0$. Moreover, $h_b=N(r,b)-r_b$, whence $$((h_{b-1}-s_{b-1})_{(b-1)})_{+1}^{+1}\geq N(r,b)-r_b.$$ If $$s_{b-1}=N(r,b-1)-((N(r,b)-r_b)_{(b)})_{-1}^{-1}+a,$$ for some $a>0$, then, by definition, we obtain $$((h_{b-1}-s_{b-1})_{(b-1)})_{+1}^{+1}=((((N(r,b)-r_b)_{(b)})_{-1}^{-1}-a)_{(b-1)})_{+1}^{+1}<N(r,b)-r_b,$$ a contradiction.\\
Conversely, let $s_1=...=s_{b-2}=0$ and $s_{b-1}\leq N(r,b-1)-((N(r,b)-r_b)_{(b)})_{-1}^{-1}$. By induction, it is easy to see that, for every $2\leq i\leq b-1$, $$((h_{i-1}-s_{i-1})_{(i-1)})_{+1}^{+1}=N(r,i).$$ Furthermore, $$((h_{b-1}-s_{b-1})_{(b-1)})_{+1}^{+1}\geq ((((N(r,b)-r_b)_{(b)})_{-1}^{-1})_{(b-1)})_{+1}^{+1}\geq N(r,b)-r_b.$$ By Remark 3.1, this is enough to show that the two vectors $H$ coincide, and the proof of the proposition is therefore complete.{\ }{\ }\qed \\
\\
Let us now take some time to consider the case in which the upper-bound $H$ of Proposition 2.3 is the same as that of Theorem A, i.e. when $s_0=s_1=...=s_{b-2}=0$ and $$s_{b-1}\leq N(r,b-1)-((N(r,b)-r_b)_{(b)})_{-1}^{-1}.$$ We want to see when we can achieve $H$ with some (generalized) compressed algebra.\\
Proposition 2.6 supplies an answer only when $s_{b-1}=0$. Instead, using sums of powers of linear forms, we can show the following:\\
\\
{\bf Theorem 3.4.} In the above hypotheses for $s$, the upper-bound $H$ of Proposition 2.3 and Theorem A is admissible (at least) for $$s_{b-1}\leq \max \lbrace N(r,b-1)-(N(r,b)-r_b),0\rbrace .$$
\\
{\bf Proof.} If $\max \lbrace N(r,b-1)-(N(r,b)-r_b),0\rbrace =0$ then we can apply Proposition 2.6. Suppose then, for the rest of the proof, that $N(r,b-1)-(N(r,b)-r_b)>0$.\\
We first show the theorem for $s_{b-1}=N(r,b-1)-(N(r,b)-r_b)$. Let us choose $s_e$ forms of degree $e$ which are the sums of $e$-th powers of $N(r,e-b)$ generic linear forms, $s_{e-1}$ forms of degree $e-1$ which are the sums of $(e-1)$-st powers of $N(r,e-1-b)$ generic linear forms, and so on for $e-2,...,b+1,b$, ending by choosing $s_b$ forms of degree $b$ which are powers of 1 generic linear form each.\\
By Theorems 2.17 and 2.18, we obtain an $h$-vector with $$h_i=\min \lbrace N(r,i),s_i\min \lbrace N(r,i),N(r,0),N(r,i-b) \rbrace $$ $$+s_{i+1}\min \lbrace N(r,i),N(r,1),N(r,i+1-b) \rbrace +...$$ $$+s_h\min \lbrace N(r,i),N(r,h-i),N(r,h-b) \rbrace +...$$ $$+s_e\min \lbrace N(r,i),N(r,e-i),N(r,e-b) \rbrace \rbrace $$ for $i=b,b+1,...,e.$\\
Since $b\leq i\leq h\leq e$, we have $h-i\leq h-b$ and, since (by Remark 2.2) $b\geq e/2$, then $h-i\leq i$. It follows that, for every $i\leq h\leq e$, $$\min \lbrace N(r,i),N(r,h-i),N(r,h-b) \rbrace =N(r,h-i).$$ Therefore, for $b\leq i\leq e$, $$h_i=\min \lbrace N(r,i),s_i+s_{i+1}N(r,1)+...+s_hN(r,h-i)+...+s_eN(r,e-i)\rbrace $$$$=\min \lbrace N(r,i),N(r,i)-r_i\rbrace =N(r,i)-r_i.$$
Notice that, for $i\geq b$, we have achieved (using sums of powers of linear forms) the same result that $[Ia]$ and $[FL]$ achieved in this range using generic forms. The advantage of using sums of powers of linear forms becomes evident when we consider the degree $b-1$. In fact, we see that, in the Inverse System, the forms of degrees greater than or equal to $b$ generate in degree $b-1$ (by Theorem 2.18) a subspace of dimension $$h_{b-1}^{'}=\min \lbrace N(r,b-1),s_b\min \lbrace N(r,b-1),N(r,1),N(r,0) \rbrace $$$$+s_{b+1}\min \lbrace N(r,b-1),N(r,2),N(r,1) \rbrace +...$$$$+s_h\min \lbrace N(r,b-1),N(r,h-b+1),N(r,h-b) \rbrace +...$$$$+s_e\min \lbrace N(r,b-1),N(r,e-b+1),N(r,e-b) \rbrace \rbrace .$$
Notice that, in each summand of the last formula, the minimum of the last two terms is always the last, i.e., for every $h$, $$\min \lbrace N(r,b-1),N(r,h-b+1),N(r,h-b) \rbrace =\min \lbrace N(r,b-1),N(r,h-b) \rbrace .$$
{\it Claim}. This minimum is always $N(r,h-b)$.\\
{\it Proof of claim}. Suppose, for some $h$ (naturally for which $s_h\neq 0$), that $N(r,b-1)<N(r,h-b)$. By hypothesis, $$N(r,b-1)>N(r,b)-r_b=s_bN(r,0)+s_{b+1}N(r,1)+...$$$$+s_hN(r,h-b)+...+s_eN(r,e-b)>s_hN(r,b-1),$$ which is a contradiction. This proves the claim.\\
Then, putting all this together, we obtain $$h_{b-1}^{'}=N(r,b)-r_b=h_b.$$
So, if $s_{b-1}=N(r,b-1)-(N(r,b)-r_b)$, we can simply choose any set of $s_{b-1}$ linearly independent forms outside the subspace described above and we are done for this case.\\
Now we show the theorem for all the other values of $s_{b-1}\leq N(r,b-1)-(N(r,b)-r_b)$. Take, for $h=b,...,e$ and $k=1,...,s_h$, non-negative integers $t_{h,k}$ such that $$N(r,h-b)+t_{h,k}\leq N(r,h-b+1)$$ and $$N(r,b-1)-(N(r,b)-r_b+\sum_{h,k}t_{h,k})>0.$$
Let us choose, for every $h=b,...,e$, exactly $s_h$ forms of degree $h$ which are, respectively, the sums of powers of $N(r,h-b)+t_{h,1}$, ..., $N(r,h-b)+t_{h,s_h}$ generic linear forms.\\
Reasoning as above, in degrees greater than or equal to $b$ we obtain again $N(r,i)-r_i$ derivatives, and in degree $b-1$ (since, similarly, we have $N(r,h-b)+t_{h,k}\leq N(r,b-1)$ for every $h$ and $k$) we obtain a subspace of dimension $$h_{b-1}^{'}+\sum_{h,k}t_{h,k}=N(r,b)-r_b+\sum_{h,k}t_{h,k}.$$ Adding $$s_{b-1}=N(r,b-1)-(N(r,b)-r_b+\sum_{h,k}t_{h,k})$$ linearly independent forms of degree $b-1$ outside the subspace described above, we are done for these values of $s_{b-1}$.\\
Notice that, in this way, we have considered all the values of $s_{b-1}\leq N(r,b-1)-(N(r,b)-r_b)$. In fact, for the least possible $s_{b-1}$ such that $r_{b-1}<0$, i.e. $$s_{b-1}=N(r,b-1)-N(r,1)s_b-N(r,2)s_{b+1}-...-N(r,e-b+1)s_e+1,$$ (naturally if positive), we need $$N(r,1)s_b+N(r,2)s_{b+1}+...+N(r,e-b+1)s_e-1$$ derivatives in degree $b-1$, and these can be obtained choosing each $t_{h,k}=N(r,h-b+1)-N(r,h-b)$ except, e.g., for the last one, $t_{b,s_b}$, which we take equal to $N(r,1)-N(r,0)-1$. For the higher values of $s_{b-1}\leq N(r,b-1)-(N(r,b)-r_b)$, of course, we just need now to decrease the values of the $t_{h,k}$; finally, when we arrive to $s_{b-1}=N(r,b-1)-(N(r,b)-r_b)$, as we have seen before, we choose all the $t_{h,k}$ equal to 0.\\
This completes the proof of the theorem.{\ }{\ }\qed \\
\\
{\bf Example 3.5.} Let $r=4$, $s=(0,0,0,0,s_4,3,0,1).$\\
It is easy to check that, for $s_4\geq 4$, we have $b=5$; the upper-bound given by Proposition 2.3 is $H=(1,4,10,20,35,13,4,1)$; moreover, for $4\leq s_4\leq 35-((13)_{(5)})_{-1}^{-1}=35-11=24$, this is also the $H$ supplied by Theorem A.\\
By Theorem 3.4, we know that $H$ is admissible at least for $4\leq s_4\leq 35-13=22$. Following the method suggested in the proof of the theorem, if $s_4=22$, we have to obtain $35-22=13$ derivatives in degree 4: a solution is to choose 1 form of degree 7 to be the sum of the powers of $N(4,7-5)=10$ generic linear forms, and 3 forms of degree 5 which are the 5-th power of 1 generic linear form each. This way, by Theorem 2.18 and Inverse Systems, we have achieved our $H$.\\
For $s_4=21$, we can choose the form of degree 7 as the sum of 11 powers and the 3 forms of degree 5 as above. Going on in this way, we can settle all the cases $4\leq s_4\leq 22$; for instance, for $s_4=4$, we choose the form of degree 7 as the sum of 20 powers, 2 forms of degree 5 as the sum of 4 powers and the third one as the sum of only 3 powers. This yields 31 derivatives, and the upper-bound $H$ is therefore achieved by adding $s_4=4$ more linearly independent forms in degree 4.\\
\\
{\bf Remark 3.6.} As far as we consider the case $s_0=s_1=...=s_{b-2}=0$ and $s_{b-1}\leq N(r,b-1)-((N(r,b)-r_b)_{(b)})_{-1}^{-1}$, we can see that, also in some other particular instances, the upper-bound $H$ is known to be admissible, even if $s_{b-1}>\max \lbrace N(r,b-1)-(N(r,b)-r_b),0\rbrace $.\\
For example, $H$ is always admissible if $b=e$. Indeed, in this case $H$ is generic up to $e-1$ and $h_e=s_e$: the fact that we can always find $s_e$ forms of degree $e$ yielding the right number of linearly independent first derivatives for each $s_{b-1}\leq N(r,b-1)-((N(r,b)-r_b)_{(b)})_{-1}^{-1}$ is shown by Cho and Iarrobino in $[CI]$, Thm. 1.4.\\ 
\\
{\bf Remark 3.7.} We will show later, however, that the upper-bound $H$ of Theorem A is not always admissible, even in some instances where it coincides with that of Proposition 2.3. Indeed, in Example 3.14, we will even see that there exist pairs $(r,s)$ for which Theorem 3.4 cannot be improved.\\
\\
Let us now come back to the general case, where we impose no restrictions on our socle-vector $s=(0,s_1,...,s_e)$.\\
\\
{\bf Definition-Remark 3.8.} Fix the pair $(r,s)$, where $r\geq $ min.emb.dim.$(s)$, and let the $h$-vector $H$ be as in Theorem A. Define $c$ as the largest integer such that $h_c$ is generic, and $t$ as the largest integer such that $$h_t=((h_{t-1}-s_{t-1})_{(t-1)})_{+1}^{+1}<N(r,t)-r_t,$$ where we set $(1_{(0)})_{+1}^{+1}=r$ and $((h_{-1}-s_{-1})_{(-1)})_{+1}^{+1}=1$, in order to avoid pathological cases.\\
Notice that we always have $0\leq t\leq e-1$ and $1\leq c\leq t+1$, since the function $r_d$ strictly increases with $d$.\\
\\
We are now ready for the second main result of this paper:\\
\\
{\bf Theorem B.} Let $(r,s)$ be as above, $r\geq $ min.emb.dim.$(s)$, and the upper-bound $H$ given by Theorem A. Then $H$ is admissible (at least) in the following cases:
\begin{itemize}
\item[  $i$)] $c=t+1$;
\item[  $ii$)] $c=t$ and $s_c\leq \max \lbrace N(r,c)-h_{c+1},0\rbrace $;
\item[  $iii$)] $c\leq t-1$ and $s_c\geq N(r,c)-c$.
\end{itemize}
{\ }\\
{\bf Proof.} If $c=t+1$, then $r_{t+1}=0$, whence $b=t+1$ and $s_1=...=s_{b-1}=0$. Therefore the upper-bound $$H=(1,r,...,h_t,h_{t+1}=h_b=N(r,b),h_{b+1}=N(r,b+1)-r_{b+1},...,s_e)$$ is achieved using generic forms, by Proposition 2.6.\\
Now let $c=t$ and suppose $\max \lbrace N(r,c)-h_{c+1},0\rbrace =0$. Then we have $s_i=0$ for every $i\leq c$, and $H$ is also given by Proposition 2.3. It is easy to see that $b=t+1$ if $r_b>0$ and $b=t$ if $r_b=0$; in either case $c\geq b-1$, and therefore $H$ is admissible by Proposition 2.6.\\
Then let $c\leq t$ and suppose therefore, if $c=t$, that $N(r,c)-h_{c+1}>0$ and $s_c\neq 0$. Suppose moreover that $s_t\leq h_t-h_{t+1}$.\\
Reasoning as in the proof of Theorem 3.4, choose $s_e$ forms of degree $e$ which are the sums of $N(r,e-t-1)$ powers of generic linear forms, $s_{e-1}$ forms of degree $e-1$ which are the sums of $N(r,e-1-t-1)$ powers, and so on down to $s_{t+1}$ forms of degree $t+1$ which are the sums of $N(r,0)=1$ power each. By Theorems 2.17 and 2.18, for $i\geq t+1$, we get $$h_i=\min \lbrace N(r,i), s_i\min \lbrace N(r,i),N(r,0),N(r,i-t-1) \rbrace $$$$+s_{i+1}\min \lbrace N(r,i),N(r,1)N(r,i+1-t-1) \rbrace +...$$$$+s_h\min \lbrace N(r,i),N(r,h-i),N(r,h-t-1) \rbrace +...$$$$+s_e\min \lbrace N(r,i),N(r,e-i),N(r,e-t-1) \rbrace \rbrace .$$
Since $b\leq t+1\leq i\leq h\leq e$, by Remark 2.2, we have $h-i\leq i$ and $h-i\leq h-t-1$. Therefore $h_i=N(r,i)-r_i$ for $i\geq t+1$, since $b\leq i$.\\
The forms of degrees higher than $t$ that we have chosen above generate, in degree $t$, a subspace of dimension $$h_t^{'}= \min \lbrace N(r,t),s_{t+1}\min \lbrace N(r,t),N(r,1),N(r,0) \rbrace +...$$$$+s_h\min \lbrace N(r,t),N(r,h-t),N(r,h-t-1) \rbrace +...$$$$+s_e\min \lbrace N(r,t),N(r,e-t),N(r,e-t-1) \rbrace \rbrace .$$
{\it Claim}. $$\min \lbrace N(r,t),N(r,h-t),N(r,h-t-1) \rbrace =N(r,h-t-1)$$ for $h=t+1,...,e$.\\
{\it Proof of claim}. Suppose that, for some $h$, $\min \lbrace N(r,t),N(r,h-t),N(r,h-t-1) \rbrace =N(r,t)$. Then $h_t^{'}=N(r,t)$, whence we have $h_t=N(r,t)$, $c=t$ and $s_c=0$, a contradiction. This proves the claim.\\
From the claim, since $t+1\geq b$, we obtain $$h_{t}^{'}=N(r,t+1)-r_{t+1}=h_{t+1}.$$  
With an argument similar to that of Theorem 3.4, suitably increasing at each step the number of summands up to $N(r,e-t)$ in degree $e$, ..., $N(r,1)$ in degree $t+1$, by Theorem 2.18 we still obtain $h_i=N(r,i)-r_i$ for $i=t+1,...,e$  and, moreover, (since we are not in the case $t=c$ and $s_c=0$) we can get any number of derivatives in degree $t$ between $h_t^{'}=h_{t+1}=N(r,t+1)-r_{t+1}$ and $$N(r,1)s_{t+1}+N(r,2)s_{t+2}+...+N(r,e-t)s_e=N(r,t)-r_t-s_t.$$ Therefore, for each $s_t\leq h_t-h_{t+1}$, we can achieve $h_t$, since $$(N(r,t)-r_t-s_t)+s_t=N(r,t)-r_t>h_t,$$ by the definition of $t$.\\
If $c=t$ there is nothing else to prove. Then, from now on, let $c\leq t-1$; suppose, moreover, that $s_t\leq h_t-h_{t+1}$, $s_{t-1}=h_{t-1}-h_t$, ..., $s_c=h_c-h_{c+1}$.\\
In order to obtain $H$, now it is enough to add $s_i$ $i$-th powers of one generic linear form each in degree $i$, for $i=c+1,...,t$. In fact, for every choice of the forms of degrees higher than $t$ that we made above in order to reach $h_t$, by Theorem 2.18, the number of derivatives yielded by those forms stabilizes in degrees less than or equal to $t$; thus, since $c<t$ and $h_c\geq h_{c+1}\geq ...\geq h_{t+1}$, by a computation similar to the one we made above, we obtain the desired values for the $h_i$.\\
Therefore the construction of a generalized compressed algebra with $h$-vector $H$ is complete under the hypotheses $s_t\leq h_t-h_{t+1}$, $s_{t-1}=h_{t-1}-h_t$, ..., $s_c=h_c-h_{c+1}$.\\
To complete the proof, now it is enough to show, for $c\leq t-1$, that $s_c\geq N(r,c)-c$, i.e. $h_c-s_c\leq c$, implies $s_t\leq h_t-h_{t+1}$, $s_{t-1}=h_{t-1}-h_t$, ..., $s_c=h_c-h_{c+1}$ (actually they will be equivalent).\\
Observe that, for $p=c,...,t-1$, $$h_{p+1}=((h_p-s_p)_{(p)})_{+1}^{+1},$$ and therefore the equality $s_p=h_p-h_{p+1}$ is equivalent to $h_p-s_p=((h_p-s_p)_{(p)})_{+1}^{+1}$, which holds if and only if $h_p-s_p\leq p$.\\
Thus it remains to show that $h_c-s_c\leq c$ implies $h_{t+1}\leq h_t-s_t$ and $h_p-s_p\leq p$, for $p=c,...,t-1$. But, if $h_c-s_c\leq c$, then $$h_{c+1}=((h_c-s_c)_{(c)})_{+1}^{+1}=h_c-s_c,$$ whence $$h_{c+1}-s_{c+1}=h_c-s_c-s_{c+1}\leq c\leq c+1.$$ By induction, we easily arrive to $$h_t=((h_{t-1}-s_{t-1})_{(t-1)})_{+1}^{+1}=h_{t-1}-s_{t-1}$$ $$=h_c-s_c-...-s_{t-1}\leq c\leq t-1< t.$$ Furthermore, $h_t-s_t\leq h_t< t$, and thus $$h_{t+1}\leq ((h_t-s_t)_{(t)})_{+1}^{+1}=h_t-s_t.$$ This completes the proof of the theorem.{\ }{\ }\qed \\
\\
{\bf Remark 3.9.} It is easy to check that the hypotheses of Proposition 2.6 are completely covered by those of Theorem B ($i)$ and $ii)$).\\
Moreover, Theorem B is a generalization of Theorem 3.4 to any socle-vector $s$. In fact, if $s_0=s_1=...=s_{b-2}=0$ and $s_{b-1}\leq N(r,b-1)-((N(r,b)-r_b)_{(b)})_{-1}^{-1}$ (i.e. the upper-bound $H$ is also given by Proposition 2.3), then it is easy to check that $t=b-1$, and $c=b-1$ if $r_b>0$ and $c=b$ if $r_b=0$. If $r_b=0$ and $c=b$, we require in both theorems that $s_{b-1}=0$, while, for $c=t=b-1$, the hypothesis $s_{b-1}\leq \max \lbrace N(r,b-1)-(N(r,b)-r_b),0\rbrace $ is of course the same as that of Theorem B, $ii$).\\
\\
{\bf Remark 3.10.} Note that, for $c\leq t-1$, $H$ must be of the following type in order to satisfy the hypotheses of Theorem B: as we saw in the proof, we must have $$h_c=N(r,c)>c\geq h_{c+1}=h_c-s_c\geq h_{c+2}=h_c-s_c-s_{c+1}\geq ...\geq $$ $$h_t=h_c-s_c-...-s_{t-1}\geq h_{t+1}\geq ...\geq h_e=s_e.$$
It is not difficult to show that the hypotheses $s_t\leq h_t-h_{t+1}$, $s_{t-1}\leq h_{t-1}-h_t$, ..., $s_c\leq h_c-h_{c+1}$, which are apparently weaker than those we worked with in the proof of Theorem B for $c\leq t-1$, i.e. $s_t\leq h_t-h_{t+1}$, $s_{t-1}=h_{t-1}-h_t$, ..., $s_c=h_c-h_{c+1}$, are actually equivalent to them. Thus, as we can see from the argument, Theorem B seems to be the best result we can show using powers of linear forms.\\
Actually, we will see in Examples 3.14 and 3.15 that, in general, the hypotheses of Theorem B, $ii$) and $iii$), cannot be improved; in fact we will exhibit pairs $(r,s)$ for which the upper-bound $H$ of Theorem A is not achieved, and such that $c=t$ and $s_c=N(r,c)+h_{c+1}+1$ in the first example and $c\leq t-1$ and $s_c=N(r,c)-c-1$ in the second.\\
\\
{\bf Example 3.11.} Let $r=3$, $s=(0,0,0,8,0,1,0,0,1)$. Then $$H=(1,3,6,10,2,2,1,1,1).$$ We have $c=3$ and $t=7$ and, since $s_3=8\geq 10-3=N(3,3)-3$, by Theorem B, $H$ is admissible. To construct a generalized compressed algebra with $h$-vector $H$, we may choose $F=L_1^8$ and $G=L_2^5$, with $L_i$ generic linear forms, and use Inverse Systems as suggested by the argument of the theorem.\\
Now let $r=3$, $s=(0,0,0,7,0,1,0,0,1).$ Then $$H=(1,3,6,10,3,3,2,2,1).$$ Since $c=3$, $t=7$ and $s_3=7\geq 10-3=N(3,3)-3$, by Theorem B, $H$ is admissible. To construct a generalized compressed algebra with $h$-vector $H$ we may choose $F=L_1^8+L_2^8$ and $G=L_3^5$, with the $L_i$ generic linear forms, and do as above.\\
\\
{\bf Proposition 3.12.} If the emb.dim. is $r=2$, then the upper-bound $H$ yielded by Theorem A is always admissible.\\
\\
{\bf Proof.} Let $r=2$, $s=(s_0=0,s_1,...,s_e).$ If $c=t+1$, then $h_t=t+1=c$ and $s_1=...=s_{c-1}=0$; moreover, $$t+2=c+1=h_c=h_{t+1}=t+2-r_{t+1},$$ whence $r_{t+1}=0$ and $b=c$; therefore we can achieve the upper-bound $H$ by Proposition 2.6.\\
If $c=t$, it is easy to see that $H$ is also given by Proposition 2.3, and thus $s_1=...=s_{b-2}=0$ and $s_{b-1}\leq b-((h_b)_{(b)})_{-1}^{-1}.$ If $h_b\leq b$, then $h_b=((h_b)_{(b)})_{-1}^{-1}$ and $H$ is achieved by Theorem 3.4. Otherwise, $h_b={b+1}$ and we are in the case $$s_{b-1}\leq b-((b+1)_{(b)})_{-1}^{-1}=b-b=0,$$ i.e. $s_{b-1}=0$, which can be settled again by Proposition 2.6, using generic forms.\\
Now let $c\leq t-1$. If $s_c\geq N(2,c)-c=c+1-c=1$, we can achieve $H$ by Theorem B. The case $s_c=0$ is clearly never verified for $c\leq t-1$, and therefore the proof of the proposition is complete.{\ }{\ }\qed \\
\\
{\bf Remark 3.13.} Proposition 3.12 may be also deduced from $[Ia]$, Thm. 4.6 C, where all the admissible $h$-vectors for a given socle-vector $s$ in emb.dim. 2 are characterized.\\
\\
As we have already mentioned, the upper-bound $H$ of Theorem A is not always admissible and, moreover, the hypotheses of Theorem B, in general, cannot be improved. We give below some examples which settle on the symmetry and the unimodality of the Gorenstein $h$-vectors in emb.dim. 3.\\
Actually, in any emb.dim. $r\geq 3$, there are examples where the $H$ of Theorem A is not admissible. Some of them, e.g., can be found thanks to the symmetry of the Gorenstein $h$-vectors; other, trivially, when $s_1>0$. In fact, in this degenerate case, in degrees greater than 1 we are actually working with $r-s_1$ variables, and therefore the admissible $h$-vectors are basically those which are admissible in emb.dim. $r-s_1$.\\ 
Furthermore, it is reasonable to believe that, as soon as something more will be known on the admissible $h$-vectors also for some other special socle-vectors (e.g. level in emb.dim. at least 3, etc.), other classes of examples of upper-bounds sharper than $H$ will probably be found as a consequence.\\
\\
{\bf Example 3.14.} Let $r=3$, $s=(0,0,0,5,0,0,1).$ Theorem A yields the upper-bound $$H=(1,3,6,10,6,3,1).$$ If it were admissible, by Inverse Systems, we would find a form $F$ of degree 6 giving 3 first derivatives, 6 second derivatives and only 5 third derivatives, to allow $s_3=5$. Hence, in emb.dim. 3, there would exist a non-unimodal Gorenstein $h$-vector $h=(1,3,6,5,6,3,1),$ which is impossible (see $[St]$, Thm. 4.2). Therefore $H$ is not admissible. Actually, now it can be easily shown that $$H^{'}=(1,3,6,10,5,3,1)$$ is the sharp upper-bound for this pair $(r,s)$.\\
Notice that here the upper-bound $H$ is given by both Theorem A and Proposition 2.3, since $b=4$, $s_1=s_2=0$ and $5=s_{b-1}\leq N(3,b-1)-((N(3,b)-r_b)_{(b)})_{-1}^{-1}=10-(6_{(4)})_{-1}^{-1}=10-5=5$. Theorem 3.4 says that this upper-bound $H$ is admissible for all the five pairs $(3,\tilde{s}=(0,0,0,\tilde{s}_{b-1},0,0,1))$ with $0\leq \tilde{s}_{b-1}\leq 10-6=4$; therefore, this example shows that Theorem 3.4 cannot be improved. Thus, not even Theorem B, $ii$) can be improved (it is easy to check that here $c=t=3$ and $s_c=5=N(r,c)-h_{c+1}+1$).\\
\\
{\bf Example 3.15.} Let $r=3$, $s=(0,0,3,0,0,0,0,1).$ Theorem A yields the upper-bound $$H=(1,3,6,4,5,6,3,1).$$ If it were admissible, by Inverse Systems, we would have a form $F$ of degree 7 with 6 second derivatives; by the symmetry of the Gorenstein $h$-vectors, $F$ should also have $6=N(3,2)$ derivatives in degree 2, which is impossible, since $s_2=3$. Therefore $H$ is not admissible. (Alternatively, reasoning as in the previous example, we can get a contradiction by observing that such an $F$ would supply a non-unimodal Gorenstein $h$-vector, which moreover is not even symmetric regardless of the value of $s_2$, since $5<6$ and $4\neq 5$).\\
Actually, we will see in $[Za]$ that this pair $(r,s)$ admits no generalized compressed algebra.\\
Notice that, in this example, $c=2$, $t=4$ (whence $c\leq t-1$) and $s_c=3=N(r,c)-c-1$. Therefore, not even Theorem B, $iii$) can be improved.\\
\\
{\bf Acknowledgements.} The author wishes to express his warm gratitude to Professor A. Iarrobino for his interesting comments on a previous version of this work.\\
\\
\\
{\bf \huge References}\\
\\
$[BG]$ {\ } A.M. Bigatti and A.V. Geramita: {\it Level Algebras, Lex Segments and Minimal Hilbert Functions}, to appear in Comm. in Algebra, 2002.\\
$[CI]$ {\ } Y.H. Cho and A. Iarrobino: {\it Hilbert Functions and Level Algebras}, J. Algebra, 241, 745-758, 2001.\\
$[EI]$ {\ } J. Emsalem and A. Iarrobino: {\it Some zero-dimensional generic singularities; finite algebras having small tangent space}, Compositio Math., 36, 145-188, 1978.\\
$[FL]$ {\ } R. Fr\"oberg and D. Laksov: {\it Compressed Algebras}, Conference on Complete Intersections in Acireale, Lecture Notes in Mathematics, 1092, 121-151, Springer-Verlag, 1984.\\
$[Ge]$ {\ } A.V. Geramita: {\it Inverse Systems of Fat Points: Waring's Problem, Secant Varieties and Veronese Varieties and Parametric Spaces of Gorenstein Ideals}, Queen's Papers in Pure and Applied Mathematics, No. 102, The Curves Seminar at Queen's, Vol. X, 3-114, 1996.\\
$[Ia]$ {\ } A. Iarrobino: {\it Compressed Algebras: Artin algebras having given socle degrees and maximal length}, Trans. Amer. Math. Soc., 285, 337-378, 1984.\\
$[Ma]$ {\ } F.H.S. Macaulay: {\it The Algebraic Theory of Modular Systems}, Cambridge Univ. Press, Cambridge, U.K., 1916.\\
$[MN]$ {\ } J. Migliore and U. Nagel: {\it Reduced arithmetically Gorenstein schemes and simplicial polytopes with maximal Betti numbers}, preprint.\\
$[St]$ {\ } R. Stanley, {\it Hilbert functions of graded algebras}, Adv. Math., 28, 57-83, 1978.\\
$[Za]$ {\ } F. Zanello: {\it Extending the idea of compressed algebra to arbitrary socle-vectors, II: cases of non-existence}, accepted for publication in J. of Algebra.

}

\end{document}